\input amsppt.sty
\magnification1200
\NoBlackBoxes
\TagsOnRight
\topmatter
\title{On Unit Fractions with Denominators in Short Intervals}\endtitle
\endtopmatter

\centerline{\smc Ernest S. Croot III}
\medskip

\centerline{\it Department of Mathematics}
\centerline{\it The University of Georgia}
\centerline{\it Athens, GA  30602}
\vskip.25in
\centerline{\it Dedicated to the memory of Paul Erd\H os}
\vskip.25in
{\eightpoint\narrower
\noindent{\bf Abstract:}
In this paper we pove that for any given rational $r>0$ and all $N > 1$,
there exist integers $N < x_1 < x_2 < \cdots < x_k < e^{r + o(1)}N$ such
that
$$
r = \frac{1}{x_1} + \frac{1}{x_2} + \cdots + \frac{1}{x_k}.
$$
}
 
\document

\heading{I.  Introduction}\endheading

Erd\H os and Graham (see \cite{3} and \cite{4}) 
asked the following questions:\bigskip

1.  Do there exist infinitely many sets of positive integers $\{x_1,x_2,...,
x_k\}$, $k$ variable, $2 \leq x_1 < x_2 < \cdots < x_k$, with
$$
1 = \frac{1}{x_1} + \frac{1}{x_2} + \cdots + \frac{1}{x_k},
$$
where $x_k/x_1$ is bounded?\smallskip

2.  If question 1 is true, what is the $\text{lim inf }x_k/x_1$ over all
such sets of integers?  Trivially, we have that this $\text{lim inf}$
is $\geq e$.  Is it actually equal to $e$?\bigskip\noindent
In this paper we will prove the following theorem, which
gives complete answers to these questions of Erd\H os and Graham. 
 
\proclaim{Main Theorem} Suppose that $r>0$ is any given rational number. 
Then, for all $N > 1$, there exist integers 
$$
N < x_1 < x_2 < \cdots < x_k \leq \left ( e^r + O_r\left ( \frac{\log\log{N}}
{\log{N}} \right ) \right )N
$$
such that 
$$
r = \frac{1}{x_1} + \frac{1}{x_2} + \cdots + \frac{1}{x_k}.
$$
Moreoever, the error term $O_r(\log\log{N}/\log{N})$ is best possible.
\endproclaim

We will now discuss the idea of the proof of the Main Theorem.  
To begin, let us suppose that
we are given some rational number $r > 0$ and an integer $N > r$.
Let $M$ be the smallest integer where
$$
r \leq \sum_{N \leq n \leq M} \frac{1}{n} \leq r + \frac{1}{M}.
$$
Using the fact that $\sum_{1 \leq n \leq t} \frac{1}{n} = \log{t} + \gamma 
+ O(1/t)$ one can show that $M=e^rN + O_r(1)$.  Now suppose
$$
\frac{u}{v} = \sum_{N \leq n \leq M} \frac{1}{n},\ \ \text{where gcd$(u,v)=1$.}
\eqno{(1)}
$$
If we had that $u/v = r$, then we would have proved our theorem for
this instance of $r$ and $N$, because $M=(e^r + O_r(1/N))N$
is well within the error of $O_r(\log\log{N}/\log{N})$ claimed by our theorem.
Unfortunately, for large $N$ it will not be the case that $u/v=r$.\smallskip

To prove the theorem, we first will use a Proposition which says that
we can remove terms from the sum in (1),
call them $1/n_1,1/n_2,...,1/n_k$, so that if
$$
\frac{u'}{v'} = \frac{u}{v} - \left \{ \frac{1}{n_1} + \frac{1}{n_2}
+ \cdots + \frac{1}{n_k} \right \} = \sum_{N \leq n \leq M \atop
n \neq n_1, n_2, ..., n_k} \frac{1}{n},\text{ where gcd$(u',v')=1$,}
$$
then all the prime power factors of $v'$ are $\leq N^{1/4 - o(1)}$,
and moreoever
$$
\frac{1}{n_1} + \frac{1}{n_2} + \cdots + \frac{1}{n_k} \asymp_r \frac{\log\log
{N}}{\log{N}}.
$$
The main idea for proving this Proposition can be found in \cite{2}, \cite{5},
and \cite{6}. 
We will then couple this with another Proposition which says that if $s$ is
some
rational number whose denominator has all its prime power factors 
$\leq N^{1/4-o(1)}$, and if $s > \frac{f(M)}{\log{M}}$, where $f(M)$ is
any function tending to infinity with $M$, then there are 
integers $M < m_1,m_2,...,m_l < e^{(c+o(1))s}M$, where $c$ is some
constant, such that
$$ 
s = \frac{1}{m_1} + \frac{1}{m_2} + \cdots + \frac{1}{m_l}.
$$
The way we use this second Proposition is we let 
$$
s = r - \frac{u'}{v'} \asymp_r \frac{\log\log{M}}{\log{M}},
$$
and then all the prime power factors of the denominator of $s$ will be
$\leq N^{1/4-o(1)} = M^{1/4-o(1)}$.  Thus, we can find
our integers $M < m_1 < \cdots < m_l < e^{(c+o(1))s}M$ as 
described above.  This will give us a unit fraction representation for
$r$ as follows:
$$
r = \sum_{N \leq n \leq M \atop n \neq n_1,n_2,...,n_k} \frac{1}{n}
\ \ +\ \  \sum_{i=1}^l \frac{1}{m_i}.
$$
All the denominators of these unit fractions will be no larger than
$$
e^{(c+o(1))s} M = e^{(c+o(1))s +r}N = \left ( e^r + O_r\left (
\frac{\log\log{N}}{\log{N}} \right ) \right ) N,
$$
and of course no smaller than $N$.  The way we will prove that the error term
$O_r\left ( \log\log{N}/\log{N}\right )$ is best-possible is by showing that
if 
$$
r = {1 \over x_1} + \cdots + {1 \over x_k},\ 2 \leq x_1 < \cdots < x_k\ 
\text{are integers,}
$$
then none of the $x_i$'s can be divisible by a prime $p>x_k/\log{x_k}$
(this idea appears in \cite{2}, \cite{3}, and  \cite{6}). 
It will turn out that this forces 
$$
{x_k \over x_1} > 
e^r\left ( 1 + {(r+o(1))\log\log{x_k} \over \log{x_k}} \right ),
$$ 
thus finishing the proof of the Main Theorem.
\smallskip

We will now state these Propositions more formally and discuss their proofs.
Before we do this, we will need the following two definitions.
Define 
$$
S(N,y)\ :=\ \{n \leq N\ :\ p^a | n \Longrightarrow 
p^a \leq y\},
$$ 
and let 
$$
\psi'(N,y) = |S(N,y)|,
$$
the number of elements in $S(N,y)$.  Our first Proposition, then, as mentioned
above is as follows:
\proclaim{Proposition 1}  Let $c > 1$ and $0 < \epsilon < \frac{1}{4}$
be given constants.  Then, for all $N$ sufficiently large, there exist integers
$$
N \leq d_1 < d_2 < \cdots < d_l \leq cN,
$$
such that if
$$
\frac{f}{g} = \sum_{N < n < cN\atop n \neq d_1,d_2,...,d_l} \frac{1}{n},
\eqno{(2)}
$$
then all the prime power factors of $g$ are $\leq N^{1/4-\epsilon}$, and
$$
\frac{1}{d_1} + \frac{1}{d_2} + \cdots + \frac{1}{d_l} = \left ( 
3\log{c} + o(1) \right ) \frac{\log\log{N}}{\log{N}}.\eqno{(3)}
$$
\endproclaim
The proof of this Proposition rests on a highly technical corollary
to a lemma taken from an earlier paper by the author (see \cite{2}).  
For completeness, we will prove both this lemma and its corollary in 
section II of the paper.
\proclaim{Lemma 1}  For all $\epsilon > 0$, there exists 
$N_\epsilon > 0$ such that if $n > N_\epsilon$ and $k > \log^
{3+2\epsilon}{n}$, then for any set of $k$ distinct primes 
$2 \leq p_1 < p_2 <\cdots <p_k < \log^{3+3\epsilon}{n}$ 
which do not divide $n$ there is a subset
$$
\{q_1,q_2,...,q_t\} \subseteq \{p_1,p_2,...,p_k\}
$$
such that 
$$
\frac{1}{q_1} + \frac{1}{q_2} + \cdots + \frac{1}{q_t} \equiv r \pmod{n},
$$ 
for any given $r$ with $0 \leq r \leq (n-1)$.
\endproclaim
\proclaim{Corollary to Lemma 1} Suppose $c > 1$, $0 < \epsilon <
\frac{1}{4}$, and $\delta > 0$ are given constants.  There exists 
a number $N_{c,\delta,\epsilon}$ so that if $N>N_{c,\delta,\epsilon}$,
then for any prime power $q$ with 
$N^{\frac{1}{4} - \epsilon} < q \leq \frac{N}{\log^{3+\delta}{N}}$ and
any residue class $r \pmod{q}$, 
there are integers $n_1,...,n_k$ satisfying: 
$$
N \leq n_1 < n_2 < \cdots < n_k \leq cN,\eqno{(4)}
$$
$$
\text{\rm where}\ \ 
n_i = qm_i,\ \text{gcd}(q,m_i)=1,\ m_i \in S(cN,q-1),\eqno{(5)}
$$
$$
\text{\rm with}\ \ 
\frac{1}{m_1} + \frac{1}{m_2} + \cdots + \frac{1}{m_k} \equiv r \pmod{q},
\eqno{(6)}
$$
$$
\text{\rm and}\ \ 
\frac{1}{n_1} + \frac{1}{n_2} + \cdots + \frac{1}{n_k} < (1+o(1))\frac{
\log^{3+2\delta/3}{N}}{N}.\eqno{(7)}
$$
\endproclaim

We will now describe how to prove our Proposition 1 using this corollary.
First, let $\delta > 0$ be some constant.  Let $d_1, d_2, ...,d_t$
be all those integers in our interval $(N,cN)$ which have a prime
power factor $q > N/\log^{3+\delta}{N}$.  Now if we let
$$
\frac{f_0}{g_0} = \sum_{N < n < cN \atop n \neq d_1,...,d_t} \frac{1}{n},
\ \ \ \text{gcd $(f_0,g_0)=1$,}
$$
then one can easily show that all the prime power factors of $g_0$ must
be $\leq N/\log^{3+\delta}{N}$.  Also, one can show that
$$
\frac{1}{d_1} + \frac{1}{d_2} + \cdots + \frac{1}{d_t} = 
\left ( (3 + \delta)\log{c}  + o(1) \right ) \frac{\log\log{N}}{\log{N}},
$$
which is a direct consequence of the following lemma:

\proclaim{Lemma 2} For $c > 1$ and $\alpha > 0$ we have 
$$
\sum_{N < mp^a \leq cN \atop p^a  > \frac{N}{\log^\alpha{N}},
\ \text{$p$ prime}} \frac{1}{mp^a} = 
\frac{\alpha(\log{c})(\log\log{N})}{\log{N}} + O\left ( \frac{1}{\log{N}}
\right ).
$$
Also, we have that
$$
\sum_{N < mp \leq cN \atop p  > \frac{N}{\log^\alpha{N}},
\ \text{$p$ prime}} \frac{1}{mp} = 
\frac{\alpha(\log{c})(\log\log{N})}{\log{N}} + O\left ( \frac{1}{\log{N}}
\right ).
$$
\endproclaim

So far we have not picked so many $d_j$'s as to violate the upper bound 
(3) claimed in Proposition 1,
since $\delta > 0$ can be chosen as small as desired; however,
the prime power factors of $g_0$ can be much larger than $N^{1/4-\epsilon}$.
Thus, the remaining numbers we choose,
$d_{t+1},...,d_l$, will have to have the properties:  if
$$
\frac{f}{g} = \frac{f_0}{g_0} - \frac{1}{d_{t+1}} - \cdots - \frac{1}{d_l}
= \sum_{N < n < cN \atop n \neq d_1,...,d_l} \frac{1}{n},\ \ \ \text{gcd
$(f,g)=1$,}
$$
then all the prime power factors of $g$ are $\leq N^{1/4-\epsilon}$, and
$$
\frac{1}{d_{t+1}} + \cdots + \frac{1}{d_l} = o \left ( \frac{\log\log{N}}
{\log{N}} \right ).\eqno{(8)}
$$

To find $d_{t+1},...,d_l$, we first select the largest prime power $q_1 | g_0$,
where $N^{1/4-\epsilon} < q_1 \leq  N/\log^{3+\delta}{N}$.  If no such prime
power exists, then we have found our integers $d_1,...,d_l$, where $l=t$,
which give rise to the property that all the prime power factors of 
$g$ are $\leq N^{1/4-\epsilon}$ (where $g$ is given by (2) above).
On the other hand, if such a $q_1=p^a$ does exist, then first write
$g_0 = q_1r_1$, where $p \nmid r_1$.  Using the Corollary to Lemma 1, 
let $d_{t+1} = n_1, d_{t+2} = n_2 , ..., d_{t+k} = n_k$, where the $n_i$'s
are as in (4) through (7) with the choices $q=q_1$ and $r=f_0/r_1$.
These new $d_j$'s are distinct from $d_1,...,d_t$, since their largest prime
factor is $q_1$, and if we let 
$$
\eqalign{
\frac{f_1}{g_1}\ &:=\ \sum_{N < n < cN \atop n \neq d_1,...,d_{t+k}}
\frac{1}{n} =  \sum_{N < n < cN \atop n \neq d_1,...,d_t} \frac{1}{n} 
- \frac{1}{n_1} - \cdots - \frac{1}{n_k}\cr
&= \frac{f_0}{g_0} - \frac{1}{n_1} - \cdots - \frac{1}{n_k} = 
\frac{1}{q_1} \left ( \frac{f_0}{r_1} - \frac{1}{m_1} - \cdots -
\frac{1}{m_k} \right ),
}
$$
where gcd$(f_1,g_1)=1$,
then all the prime power factors of $g_1$ are $\leq q_1 - 1$.  To see
this, we have from (6) that if
$$
\frac{w_1}{w_2} = \frac{f_0}{r_1} - \frac{1}{m_1} - \cdots - \frac{1}{m_k},
\ \ \text{ gcd$(w_1,w_2)=1$,}
$$  
then $p | w_1$, and so $q_1 \nmid g_1$ (and the same goes for any prime power
bigger than $q_1$).   
From (7) we have one final property that our $d_j$'s satisfy: 
$$
{1 \over d_{t+1}} + {1 \over d_{t+2}} + \cdots + {1 \over d_{t+k}}
= {1 \over n_1} + \cdots + {1 \over n_k} < 
(1 + o(1)){\log^{3+2\delta/3}{N} \over N}.
$$

We now repeat the process as above and select the largest prime power factor of
$g_1$, call it $q_2$, where $N^{1/4-\epsilon} < q_2 \leq q_1-1$.  If no such
prime power exists, then we are finished and have found our integers
$d_1,...,d_l$ with $l=t+k$.  If such $q_2$ does exit, we can use lemma 2
again as we did above to find our integers $d_{t+k+1},d_{t+k+2},...,
d_{t+k+h}$ in $(N,cN)$, distinct from $d_1,...,d_{t+k}$, such that if we let
$$
{f_2 \over g_2}\ :=\ \sum_{N < n < cN\atop n\neq d_1,...,d_{t+k+h}}
{1 \over n},
$$
where gcd$(f_2,g_2)=1$, then the largest prime power factor of $g_2$ is
at most $q_2-1$.  Also, 
$$
{1 \over d_{t+k+1}} + \cdots + {1 \over d_{t+k+h}} < (1+o(1)){\log^{3+2\delta/3}
{N} \over N}.
$$

If we continue in this manner of picking $d_j$'s to cancel off prime power 
factors $> N^{1/4-\epsilon}$, we will eventually find our integers 
$d_1,...,d_l$ such that if $f$ and $g$ are as in (2), then all the
prime power factors of $g$ are $\leq N^{1/4-\epsilon}$.  To see that 
our $d_j$'s satisfy (8), and therefore (3), we observe that
$$
\eqalign{
\sum_{t+1 \leq j \leq l} {1 \over d_j} &< (1+o(1)){\log^{3+2\delta/3}{N} \over 
N} \sum_{N^{1/4-\epsilon} < p^a \leq N/\log^{3+\delta}{N} \atop p\ prime}
1 \cr
&=(1+o(1)){\log^{3+2\delta/3}{N} \over N} \pi(N/\log^{3+\delta}{N})
= {1 + o(1)\over \log^{1+\delta/3}{N}}.
}
$$

We now formally state our second Proposition mentioned above and describe
its proof:  
\proclaim{Proposition 2} Suppose $a,\ b$ are positive integers, 
where gcd$(a,b)=1$, all the prime power factors of $b$ are $\leq M^{1/4-
\epsilon}$, where $0 < \epsilon < 1/8$.  Futher, we will allow the
size of $a/b$ to depend on $M$:  
suppose ${f(M) \over \log{M}} < \frac{a}{b} \leq 1$, where 
$f(M)<\log{M}$ is any function tending to infinity with $M$.  
Select $c(M) > 0$ such that
$$
2\frac{a}{b} \leq \sum_{M \leq n \leq c(M)M \atop
n \in S(M,M^{\frac{1}{4} - \epsilon})} \frac{1}{n} < 2\frac{a}{b} + {1 \over
c(M)M}.
$$
Remark:  We will show that $c(M) = e^{(v(\epsilon) + o(1))a/b}$, where
$v(\epsilon)$ is some function depending only on $\epsilon$. 
Then for all $M$ sufficiently large, there exist integers
$$
M \leq n_1 < n_2 < \cdots < n_k \leq c(M)M,
$$
each $n_i \in S(M,M^{1/4-\epsilon})$ such that 
$$
\frac{a}{b} = \frac{1}{n_1} + \frac{1}{n_2} + \cdots + \frac{1}{n_k}.
$$
\endproclaim

Let $m_1,...,m_l$ be all the integers where
$$
M \leq m_1 < m_2 < \cdots < m_l \leq c(M)M,\ m_j \in S(c(M)M,M^{1/4-\epsilon}).
$$
It will turn out that
$$
l \gg_{a,b,\epsilon} M.
$$
The proof of Proposition 2 rests entirely on estimating the following
exponential sum:
$$
E\ :=\ \sum_{h=-P/2}^{P/2-1} e(-ah/b)A(h),
$$
where $e(\cdot) := e^{2\pi i \cdot}$,
$$
A(h)\ :=\ \prod_{j=1}^l \left ( 1 + e(h/m_j) \right ),
$$
and 
$$
P\ :=\ \text{lcm}\{2,3,4,...,[N^{1/4-\epsilon}]\}.
$$
It turns out that 
$$
\#\{ \{n_1,...,n_k\} \subseteq \{m_1,...,m_l\}, k\ \text{variable}\  
:\ 1/n_1 + \cdots + 1/n_k = a/b\} \geq {E \over P} - 2.
$$
The $-2$ comes from the fact that in the case $a/b=1$, the exponential sum
picks up the extraneous representations for $a/b=0$ and $a/b=2$, and there
can be at most one such representation each.  In the cases where $a/b<1$,
we can omit the $-2$ above to get the exact count:
$$
\#\{ \{n_1,...,n_k\} \subseteq \{m_1,...,m_l\}, k\ \text{variable}\  
:\ 1/n_1 + \cdots + 1/n_k = a/b\} = {E \over P}.
$$

The way we obtain a lower bound for the exponential sum $E$ is by showing:
\bigskip\indent
1.  For $|h| < M/2$, $\text{Re}(e(-ah/b)A(h)) > 0$; and so,
$$
\text{Re}\left \{\sum_{|h| < M/2} e(-ah/b)A(h) \right \} =
A(0) + \text{Re}\left \{\sum_{|h| < M/2,\  h \neq 0} e(-ah/b)A(h)\right \}
\geq A(0) = 2^l.
$$

2.  For $|h| \geq M/2$ and $|h| \leq P/2$, $|A(h)| < {2^{l-1} \over P}$; 
and so,
$$
\sum_{|h| \geq M/2,\ |h| \leq P/2} |A(h)| < 2^{l-1}.
$$
\bigskip\noindent
Putting together these two facts, we find that our number of representations
for $a/b$ is at least 
$$
{|E| \over P} - 2 \geq {2^{l-1} \over P} - 2 > 0,
$$
since
$$
2^l \gg_{a,b,\epsilon} 2^{cM},
$$
for some constant $c$, while
$$
P < e^{M^{1/4-\epsilon}(1+o(1))}. 
$$

\heading{II.  Technical Lemmas and Their Proofs}\endheading

\demo{Proof of Lemma $1$}
Suppose that $b$ is coprime to $n$ and let 
$r_n(a/b)$ denote the least residue of $ab^{-1} \pmod{n}$ in absolute value.    
The number of subsets of $\{ p_1, ..., p_k\}$ whose
sum of reciprocals is $\equiv l \pmod{n}$ is then given by
$$
S_l\ :=\ \frac{1}{n} \sum_{h=0}^{n-1} e\left ( \frac{-hl}{n}\right )
\prod_{j=1}^k \left ( 1+ e\left(\frac{r_n(h/p_j)}{n}
\right ) \right ),
$$
where $e(x)$ is defined to be $e^{2\pi ix}$.  Define
$$
P(h)\ :=\ \prod_{j=1}^k \left ( 1 + e\left ( \frac{r_n(h/p_j)}{n} \right )
\right ).
$$
We will show that
$$
| P(h) | < \frac{2^k}{n}, \eqno{(8)}
$$
when $h \neq 0$ and when $n$ is sufficiently large.  It will then follow that
$$
|S_l| = \left | \frac{1}{n} \sum_{h=0}^{n-1} P(h) \right |
> \frac{1}{n} \left\{ 2^k - \sum_{h=1}^{n-1} \frac{2^k}{n} \right \}
= \frac{2^k}{n^2} > 0,
$$
and thus there is at least one subset of $\{p_1,...,p_k\}$ with the
desired property.\smallskip 
To prove (8) we note that 
$$
\eqalign{
| P(h)|  &= \left |\prod_{j=1}^k \left ( 1+ e\left(\frac{r_n(h/p_j)}
{n} \right ) \right )\right |\cr
&= \left |\prod_{j=1}^k \left ( e\left(-\frac{r_n(h/p_j)}{2n}\right )
+ e\left(\frac{r_n(h/p_j)}{2n} \right ) \right )\right |\cr 
&= 2^k\prod_{j=1}^k \left |\cos\left (\pi\frac{r_n(h/p_j)}{n} \right )
\right |}\eqno{(9)}
$$

We may write
$$
r_n ( h/p_j) = \frac{s_jn + h}{p_j},
$$
where $0 \leq h \leq (n-1)$ and $s_j$ is an integer satisfying 
$-\left [\frac{p_j}{2}\right ] < s_j
\leq \left [ \frac{p_j}{2} \right ]$.  Define $L(x)\ :=\ 
\log^{2+2\epsilon}{x}+1$.
We will now show that when $n$ is sufficiently large at least 
$\frac{k}{2}$ of the $s_j$'s have the property that 
$|s_j| > L(n)$:  for if we suppose there are infinitely many $n$ where at
least $\frac{k}{2}$ of the $s_j$'s satisfy $|s_j| \leq L(n)$ then, by the
pigeonhole principle, there is a number $m$ with $|m| \leq L(n)$ such that
$s_j = m$ for at least
$$
\frac{k/2}{2L(n)+1} > \frac{\log^{3+2\epsilon}{n}}{4\log^{2+2\epsilon}{n}+6}
\gg \log{n}
$$
of the primes $p_j$ dividing $mn+h$ when $n$ is sufficiently large.  However,
this is impossible for large $n$ since $|mn+h| < |n(L(n)+1)| < n^2$ has
$o(\log{n})$ distinct prime factors.  Thus when $n$ is sufficiently large
at least $\frac{k}{2}$ of the $s_j$'s satisfy $|s_j| > L(n)$.

It follows that, when $n$ is sufficiently large, at least $\frac{k}{2}$ of the
$p_j$'s satisfy
$$
| r_n(h/p_j) | = \left | \frac{s_jn+h}{p_j}
\right | > \left | \frac{(s_j-1)n}{p_j} \right | > 
\frac{n}{\log^{1+\epsilon}{n}}.
$$
We have for such primes $p_j$ that when $n$ is sufficiently large,
$$
\eqalign{
\left |\cos\left ( \pi {r_n(h/p_j) \over n} \right ) \right |
&= \left | 1 - {\pi^2 \over 2} \left ( {r_n(h/p_j) \over n} \right )^2
+ O\left ( \left ( {r_n(h/p_j) \over n} \right )^4 \right ) \right |\cr
&< 1 - {\pi^2 \over 2\log^{2+2\epsilon}n} + O \left ( {1 \over 
\log^{4+4\epsilon}n} \right ). 
}
$$
and so, from (9), since $k > \log^{3+2\epsilon}{n}$ we have that 

$$
|P(h)| < 2^k \left ( 1 - \frac{\pi^2}{\log^{2+2\epsilon}{n}} +
O\left ( \frac{1}{\log^{4+4\epsilon}{n}} \right ) \right )^{k/4}
\ll 2^ke^{-\frac{\pi^2 \log{n}}{4}} = o\left ( \frac{2^k}{n} \right ),
$$
which was just what we needed to show in order to prove our lemma. 
\enddemo
\demo{Proof of Corollary}
Let $s(q)$ denote the smallest integer with   
$$
s(q) > {N \over q\log^{3+\delta}q},\ s(q) \in 
S(cN,N^{1/4-\epsilon}),\ \text{and gcd}(q,s(q))=1.
$$
This number $s(q) = (1 + o(1)) \frac{N}{q \log^{3+\delta}{q}}$.  We will
construct the $m_i$'s so that $m_i = s(q)r_i$, where $r_i$ is a small prime.
Let
$$
\frac{N}{qs(q)} < p_1 < p_2 < \cdots < p_l < \frac{cN}{qs(q)} < (c+o(1))
\log^{3+\delta}{q}
$$ 
be all the primes
between $\frac{N}{qs(q)}$ and $\frac{cN}{qs(q)}$ which do not divide $q$.
The number of these primes is at least
$$
\eqalign{
\pi\left ( {cN \over qs(q)} \right ) - \pi\left ( {N \over qs(q)} \right ) - 1
&= \pi\left \{ (c+o(1))\log^{3+\delta}{q} \right \} - \pi\left \{
(1+o(1))\log^{3+\delta}{q} \right \} \cr
&= \left ( {c-1-o(1) \over 3+\delta}\right ) 
{\log^{3+\delta}{q} \over \log\log{q}}.
}
$$

When $N$ is sufficiently large we have from our lemma 1 above with
$\epsilon = \delta/3$ that
there is a subset $r_1 < r_2 < \cdots < r_k$ of the primes
$\{p_1,p_2,...,p_l\}$ with
$$
\frac{1}{s(q)r_1} + \cdots + \frac{1}{s(q)r_k} \equiv r \pmod{q}.
$$
where $N < qs(q)r_i < cN$ for all $i=1,2,...,k$; moreover, there is such 
a subset with $k < (1+o(1))\log^{3+\frac{2}{3}\delta}{N}$.  Thus, if we let 
$m_i = s(q)r_i$ and therefore $n_i = qm_i = qs(q)r_i$, we satisfy
(4), (5), and (6).  If we assume $k < (1+o(1))\log^{3+\frac{2}{3}\delta}{N}$,
as we are allowed to do, then 
$$
\frac{1}{n_1} + \frac{1}{n_2} + \cdots + \frac{1}{n_k} < (1+o(1)) 
\frac{\log^{3+\frac{2}{3}\delta}{N}}{cN},
$$ 
which satisfies (7).
\enddemo 

\demo{Proof of Lemma 2}  Using the the fact that $\sum_{1 \leq j \leq n}
\frac{1}{j} = \log{n} + \gamma + O(1/n)$, together with the estimate
$$
\sum_{p^a \leq n \atop \text{$p$ prime}} \frac{1}{p^a} = \log\log{n} + B +
o(1/\log{n}),
$$
where $B$ is some constant, we have the following chain of inequalities:
$$
\eqalign{
\sum_{N < mp^a \leq cN \atop p^a > \frac{N}{\log^\alpha{N}},\ \text{$p$ prime}} 
\frac{1}{mp^a} &= \sum_{{N \over \log^\alpha{N}} < p^a
\leq cN} \frac{1}{p^a} \sum_{N/p^a < m \leq cN/p^a} \frac{1}{m}\cr
&= \sum_{{N \over \log^\alpha{N}} <p^a \leq cN} \frac{1}{p^a} \left \{
\log(cN/p^a) - \log(N/p^a) + O(p^a/cN)  \right \}\cr
&= \sum_{{N \over \log^\alpha{N}} < p^a \leq cN} \frac{1}{p^a} \left \{
\log{c} + O(p^a/cN) \right \}\cr
&= \log{c} \sum_{{N \over \log^\alpha{N}} < p^a \leq cN} \frac{1}{p^a}
+ O\left ( \frac{\pi(cN)}{cN} \right ) \cr
&= \log{c} \left \{ \log\log{cN} - \log\log\left ( {N \over \log^\alpha{N}}
\right )  + o (1/\log{N}) \right \}\cr 
&\ \ \ \ \ \ \ \ \ \ \ \ \ \ \ \ \ \ \ \ + O(1/\log{N})\cr
&= \frac{\alpha(\log{c})(\log\log{N})}{\log{N}} + O(1/\log{N}),
}
$$
as claimed.  The proof for the sum over primes $p$, instead of prime
powers $p^a$, is exactly the same. 
\enddemo

\heading{III.  Proof of Proposition 1}\endheading

Fix a $\delta > 0$ and let $N^{\frac{1}{4}-\epsilon} \leq q_1 < q_2 <
\cdots < q_h < \frac{N}{\log^{3+\delta}{N}}$ be all the prime powers
between $N^{\frac{1}{4}-\epsilon}$ and $\frac{N}{\log^{3+\delta}{N}}$.
Define
$$
S\ :=\ \{N\leq n\leq cN\},
$$
$$
S_{h+1}\ :=\ S \setminus \{n\ :\ n=mp^a,\text{ where
$\frac{N}{\log^{3+\delta}{N}}\leq p^a \leq N$, $p$ prime}\},
$$
and let 
$$
\frac{u_{h+1}}{v_{h+1}} = \sum_{n \in S_{h+1}} \frac{1}{n},
$$
where gcd$(u_{h+1},v_{h+1})=1$.  We observe that all of the prime power
factors of $v_{h+1}$ are smaller than $\frac{N}{\log^{3+\delta}{N}}$ and 
by lemma 2 we have
$$
\frac{u_{h+1}}{v_{h+1}} = \sum_{N\leq mp^a \leq cN\atop p^a \geq \frac{N}
{\log^{3+\delta}{N}}} \frac{1}{mp^a} = \left ((3+\delta)\log{c}+ o(1) \right )
\frac{\log\log{N}}{\log{N}}.
$$
Starting with the prime power $q_h$ we will successively construct sets
$$
S_{h} \supseteq S_{h-1} \supseteq S_{h-2} \supseteq \cdots \supseteq S_1,
$$
where if 
$$
\frac{u_i}{v_i} = \sum_{n \in S_i} \frac{1}{n},
$$
gcd$(u_i,v_i)=1$, then all the prime power factors of $v_i$ are smaller
$q_i$, for all $i=1,2,..,h+1$; moreover, we will construct these sets in such
a way that
$$
\sum_{n \in S\setminus S_1} \frac{1}{n} = \left ( (3+\delta)\log{c} + o(1) 
\right) \frac{\log\log{N}}{\log{N}}.
$$
If we can accomplish this, then we can just let $\{d_1,...,d_l\} =
S \setminus S_1$ and satisfy the requirements of the Proposition.\smallskip

Suppose, for proof by induction, we have constructed the sets 
$S_i$ where $2 \leq i \leq h+1$.  If $q_{i-1} \nmid v_i$, we just let 
$S_{i-1}\ :=\ S_i$, and then all the prime power factors of $v_{i-1}$
are smaller than $q_{i-1}$.  On the other hand, if $q_{i-1} \nmid v_i$, then
using the corollary to lemma 1 we can find integers 
$N < n_1 < n_2 < \cdots < n_k < cN$ where $n_j = q_{i-1}m_j$, 
gcd$(q_{i-1}, m_j)=1$, all the prime power factors of the $m_j$'s are
smaller than $q_{i-1}$, 
and
$$
\frac{1}{m_1} + \cdots + \frac{1}{m_k} \equiv q_{i-1}\sum_{n \in S_i} \frac{1}{n}
= q_{i-1}\frac{u_i}{v_i} \pmod{q_{i-1}}.
$$
Then if we let $S_{i-1}\ :=\ S_i \setminus \{n_1,n_2,...,n_k\}$ we will have
that 
$$
q_{i-1}\frac{u_{i-1}}{v_{i-1}} = q_{i-1}\frac{u_i}{v_i} - 
\frac{1}{m_1} - \cdots - \frac{1}{m_k} \equiv 0 \pmod{q_{i-1}},
$$ 
and so $q_{i-1}$ does not divide $v_{i-1}$, nor does any other prime
power bigger than $q_{i-1}$ since all the prime power factors of $v_i$ and
the $n_j$'s are at most $q_{i-1}$.  We conclude, by induction, that
$S_i$ can be constructed for $1 \leq i \leq h+1$.\smallskip

From the corollary to lemma 1, for each $2 \leq i \leq h+1$ we can
pick the $n_j$'s as above so that 
$$
\sum_{n \in S_i \setminus S_{i-1}} \frac{1}{n} < (1+o(1)) \frac{\log^{3+
\frac{2}{3}\delta}{N}}{N}.
$$
It follows that
$$
\sum_{n \in S_{h+1} \setminus S_1} \frac{1}{n} < (1+o(1)) \frac{\pi\left (
\frac{N}{\log^{3+\delta}{N}} \right )\log^{3+\frac{2}{3}\delta}{N}}{N}
= \left (1+o(1)\right )\frac{1}{\log{N}},
$$
and so if we let
$$
\{d_1,d_2,...,d_l\} = S\setminus S_1,
$$
then (2) and (3) are satisfied and
$$
\frac{1}{d_1} + \cdots + \frac{1}{d_l} = 
\sum_{n \in S \setminus S_1} \frac{1}{n} = 
\left ((3+\delta)\log{c} + o(1) \right ) \frac{\log\log{N}}{\log{N}}.
$$
Since we can choose $\delta$ as small as desired, the Proposition follows.

\heading{IV.  Proof of Proposition 2}\endheading

First we will show that
$$
c(M) = e^{(v(\epsilon)+o(1))a/b}.
$$
where $v(\epsilon)$ is some constant depending only on $\epsilon$.  To
do this we will need the following lemma:

\proclaim{Lemma 3 (N.G. de Bruijn)}  For any fixed $\epsilon < 3/5$,
uniformly in the range
$$
y \geq 2,\ \ 1 \leq u \leq \exp\{(\log{y})^{3/5-\epsilon}\},
$$
we have
$$
\psi(x,y) = x\rho(u) \left \{1 + O\left ( {\log(u+1) \over 
\log{y} } \right ) \right \},
$$
where $u=\log{x}/\log{y}$ and $\rho(u)$ is the unique continuous solution
to the differential-difference equation 
$$
\cases &u\rho'(u) = -\rho(u-1),\ \ \text{if $u > 1$} \\
       &\rho(u)=1,\ \ \text{if $1 \leq u \leq 1$}.
\endcases
$$
\endproclaim
\noindent
(For a proof of this lemma, see \cite{1}.) \smallskip

Using lemma 3 with 
$$
u = {1 \over {1/4 - \epsilon}},\ \text{\rm and}\  x=M
$$
gives us that
$$
\psi(M+z,M^{1/4-\epsilon}) - \psi(M,M^{1/4-\epsilon}) \sim z\rho(u) 
$$
for $z \gg M/\log{M} $.  Using this and partial summation it is fairly easy
to see that for $c'(M) = e^{(2/\rho(u) +o(1))a/b}$,
$$
\sum_{M \leq n \leq c'(M)M \atop p | n \Longrightarrow p < M^{1/4-\epsilon}} {1
\over n} \sim 2{a \over b},
$$
for $f(M)/\log{M} < a/b < 1$, where $f(M)$ is any function tending to infinity
with $M$.  The error incurred by replacing the condition `$p|n \Longrightarrow
p < M^{1/4-\epsilon}$' with `$n \in S(c'(M)M,M^{1/4-\epsilon})$' will be at
most
$$
\eqalign{
\sum_{n \leq c'(M)M \atop {p^a | n, p^a > M^{1/4-\epsilon} \atop 
\text{where $p < M^{1/4-\epsilon}$ is prime}}} {1 \over n}
&\ll \sum_{p \leq M^{1/8-\epsilon/2} \atop p\ \text{prime}} 
{1 \over M^{1/4-\epsilon}} \sum_{m \leq c'(M)M^{3/4+\epsilon}} {1 \over m} 
\cr
&\ \ \ \ \ +\ \sum_{M^{1/8-\epsilon/2} < p \leq M^{1/4-\epsilon} 
\atop p\ \text{prime}} {1 \over p^2} \sum_{m \leq c'(M)M/p^2} {1 \over m}\cr
&\ll M^{1/8-\epsilon/2}. 
}
$$
Thus, we see that
$$
\sum_{M \leq n \leq c'(M)M \atop n \in S(c'(M)M,M^{1/4-\epsilon})} {1 \over n}
\sim 2{a \over b},
$$
which gives us that
$$
c(M) \sim c'(M) = e^{(2/\rho(u) + o(1))a/b}.
$$

Let
$$
P\ :=\ \text{ lcm}(1,2,3,...,[M^{1/4-\epsilon}] ) = \prod_{p \leq 
M^{1/4-\epsilon} \atop p\ \text{prime}} p^{a_p} = e^{M^{1/4-\epsilon}(1+o(1))}, 
$$
where $a_p$ is the largest integer such that $p^{a_p} \leq M^{1/4-\epsilon}$. 
Let $M \leq m_1 < m_2 < \cdots < m_l \leq c(M)M$ be all the divisors of $P$
lying in
$[M,c(M)M]$; that is, all the integers in $S(c(M)M,M^{1/4-\epsilon})$ in the
interval $[M,c(M)M]$. 
By standard methods of exponential sums, one has that
$$
\eqalign{
\#\{ \{n_1,...,n_k\} &\subseteq \{m_1,...,m_l\}, k\ \text{variable}\  
:\ 1/n_1 + \cdots + 1/n_k = a/b\}\cr
&\geq {1 \over P} \sum_{h=-P/2}^{P/2-1} e\left ({-ah \over b} \right ) 
\prod_{j=1}^l \left \{ 1 + e\left ( {h \over m_j} \right ) \right \}
- 2,
}
$$
where $e(\cdot) = e^{2\pi i \cdot}$. 
The reason for subtracting 2 in the above equation is that when
$a/b=1$, the exponential sum not only counts subsets summing to 1, but also
0 and 2. 

Let
$$
A(h)\ :=\ \prod_{j=1}^{l}  \left \{ 1 + e\left ( {h \over m_j} 
\right ) \right \}
= e\left ( {h \over 2} \left \{ {1 \over m_1} + \cdots + {1 \over m_l}
\right \} \right ) \left ( 2^l \prod_{j=1}^l \cos ( \pi h/m_j) \right ).
\eqno{(9)}
$$
Upon substituting in our equation above this gives   
$$
\eqalign{
\#\{\{n_1,...,n_k\} &\subseteq \{m_1,...,m_l\}
, k\ \text{variable}\  :\ 1/n_1 + \cdots + 1/n_k
= a/b\}\cr
 &\geq {1 \over P} \left (\sum_{h=-P/2}^{P/2-1} e(-ah/b) A(h)
\right ) - 2.
}\eqno{(10)}
$$ 

We will now try to find a lower bound for (10).
To do this we will show that  
$$
|A(h)| < {2^l \over 2P},\ \text{for $-P/2 \leq h \leq P/2-1$ with $
|h| > M/2$.}\eqno{(11)}
$$
and that
$$
\text{Re}\left ( \sum_{|h| \leq M/2} e(-ah/b)A(h) \right ) 
> 2^l, \eqno{(12)}
$$
From (10), (11), and (12) it then follows that
$$
\eqalign{
\#\{ \{n_1,...,n_k\} &\subseteq \{m_1,...,m_l\},\ \text{$k$ variable}\ :\ 
1/n_1 + \cdots + 1/n_k = a/b\}\cr
& > {2^{l-1} \over P} - 2 = 2^{l - O(M^{1/4-\epsilon})},
}
$$
which is exponential in $l$ since 
$$
l \gg_\epsilon M{a \over b} \gg {M \over \log{M}}.
$$ 

To establish (12), we first observe from (9) that 
$$
\text{Arg}\{e(-ah/b)A(h)\} = 
{-2\pi a h\over b} + {\pi h}\left \{ {1 \over m_1} + \cdots + {1 \over m_l}
\right \} + \text{Arg}\left \{\prod_{j=1}^l \cos(\pi h/m_j) \right \}.
\eqno{(13)}
$$
Using the fact that
$$
{1 \over m_1} + \cdots + {1 \over m_l} = 2{a\over b} + \delta,
$$
where 
$$
0 \leq \delta \leq {1 \over c(M)M},
$$
together with the fact that each $m_j$ is $\geq M$, we have 
$$
\left | {-2\pi a h \over b} + {\pi h}\left \{ {1 \over m_1} + \cdots
 + {1 \over m_l} \right \}\right | = \pi \delta |h| < {\pi |h| \over M} <
{\pi \over 2}, \eqno{(14)}
$$
whenever 
$$
|h| < {M \over 2}.
$$
Also for such $h$, we observe that
$$
\cos(\pi h/m_j) \geq \cos(\pi/2) = 0,\ \text{for $j=1,2,...,l$},
$$
since the $m_j's$ are all $\geq M$.  Using this, together with (13) and (14),
we find that
$$
|\text{Arg}\{e(-ah/b)A(h)\}| \leq {\pi \over 2},\ \text{whenever $|h| < 
{M \over 2}$.}
$$
Thus, for such $h$ we have
$$
Re\{ e(-ah/b)A(h) \} \geq 0,
$$
and so
$$
\text{Re}\left ( \sum_{|h| \leq M/2} e(-ah/b)A(h) \right )
= 2^l + \text{Re}\left ( \sum_{|h| \leq M/2 \atop h \neq 0} e(-ah/b)A(h)
\right ) \geq 2^l,
$$
which establishes (12).\indent

In order to establish (11), we will need the following lemma, which will
be proved in the next section of the paper: 

\proclaim{Lemma 4}  Suppose $0 < \epsilon < \frac{1}{8}$.
Let $M \leq m_1 < m_2 < \cdots < m_l \leq \left ( 1 + {1 \over \log{M}} 
\right )M$ be all the integers in this interval with 
$m_i \in S(M,M^{\frac{1}{4} - \epsilon})$.  Then
for $M$ sufficiently large and $h$ real, either\bigskip
\indent 1.  There are $\gg M^{\frac{3}{4}}$ $m_i$'s which do not divide any
integer in $I\ :=\ (h-M^{\frac{3}{4}},h+M^{\frac{3}{4}})$, or
\smallskip
\indent 2.  There is an integer in this interval which is divisible by
$P := lcm \{ p^a \leq M^{\frac{1}{4} - \epsilon}\ :\ p\ prime\}$.
\endproclaim

From this lemma, it follows that if
$$
{M \over 2} \leq |h| \leq P/2,
$$
then for some constants $c_1,c_2>0$ there are $> c_1M^{3/4}$ $m_j$'s 
such that for any integer $z$
$$
\left | {h \over m_j} - z \right | >{c_2 \over M^{1/4}},
$$
for all $M$ sufficiently large.
For these integers $m_j$, we will have that
$$
\left | \cos\left ( \pi h/m_j \right ) \right | < \left | \cos\left (
\pi c_2 / M^{1/4} \right ) \right | =  1\ -\ {1 \over 2} {\pi^2 c_2^2
\over M^{1/2}} \ +\ O\left ( {1 \over M} \right ).
$$
From this and (9) it follows that for such $h$ 
$$
|A(h)| < 2^l \left ( 1 - {1 \over 2} {\pi^2 c_2^2 \over M^{1/2}} +
O \left ( {1 \over M} \right ) \right )^{c_1 M^{3/4}}
\ll 2^l e^{\pi^2 c_1c_2^2 M^{1/4}/2} = o\left ( {2^l \over P} \right ).
$$
This establishes (11) and thus proves the Proposition.

\heading{V.  Proof of Lemma 4}\endheading
 
For each integer $n$ satisfying 
$$
M^{\frac{3}{4}}\log^2{M} < n < 2M^{\frac{3}{4}}\log^2{M},
\ \ \text{\rm and }\ n \in S(2M^{\frac{3}{4}}\log^2{M},
M^{1/4 - \epsilon}), \eqno{(15)} 
$$
define 
$$
M(n)\ :=\ \{ m_j\ :\ m_j = nq,\ \omega(q) \leq 3 \}.
$$
We claim that $\text{\rm lcm }\ M(n)= P$ for all such $n$.  
We will show below that the truth 
of this claim implies that either:  \bigskip

A.  There is an $n$ satisfying (15) such that every integer of
$M(n)$ divides a single integer in $I$, which
together with the assumption $\text{ lcm } M(n)=P$, gives us case 2 in
the claim of our lemma, or\smallskip

B.  For each $n$ satisfying (15), there is an integer 
$m_{\alpha(n)} \in M(n)$ which does not divide any integer in
$(h-M^{3/4},h+M^{3/4})$.\bigskip\noindent
We will assume that case B is true and show that it implies case 1
in the claim of our lemma (and thus if we can show that $\text{lcm } M(n)=P$
and that either A or B is true, we may conclude 
that either case 1 or case 2 in our lemma is true):\smallskip

The first thing to notice is that from Lemma 3 we know there are
$\gg_\epsilon M^{3/4}\log^2{M}$ integers $n$ satisfying (15).   
If all of the $m_{\alpha(n)}$'s as indicated in case B were distinct,
then we would have that there are $\gg_\epsilon M^{3/4}\log^2{M}$ 
$m_j$'s not dividing any integer in $(h-M^{3/4},h+M^{3/4})$, which is
the first possibility claimed by our lemma;  however, it is not necessarily
the case that the $m_{\alpha(n)}$'s are distinct.  To overcome this 
difficulty, we will now show that no $m_i$ can live in too many of the
sets $M(n)$:  Let
$$
\eqalign{
D(M)\ &:=\ \max_{m_i}\ \#\{n\ :\ \text{$n$ satisfies (15) and $m_i \in M(n)$}\}\cr
&\leq \max_{m_i}\ \#\{q\ :\  q | m_i, \omega(q) \leq 3, 
q\geq \frac{M^{1/4}}{2\log^2{M}} \} = o\left ( \log^2{M} \right ),
}
$$
then
$$
\#\{m_{\alpha(n)}\ :\ \text{$n$ satisfies (15)} \} \geq
\frac{\psi(2M^{\frac{3}{4}}\log^2{M},M^{\frac{1}{4} - \epsilon})
- \psi(M^{\frac{3}{4}}\log^2{M}, M^{\frac{1}{4}-\epsilon})}
{D(M)} \gg M^{\frac{3}{4}}. 
$$
Thus, there are $\gg M^{\frac{3}{4}}$
$m_j$'s which do not divide any integer in $(h-M^{3/4},h+M^{3/4})$, which
covers case 1 claimed by our lemma.\smallskip

We now will show that if $\text{lcm }M(n) = P$ for all $n$ satisfying (15),
then either case A or case B above must be true.  So, let us assume then that
$\text{lcm }M(n) = P$ for all $n$ satisfying (15).  If case B is true, then
we are done.  So, let us assume that case B is false.  Then, we must have
there there is an $n$ satisfying (15) such that each member of $M(n)$ divides
an integer in $I$.  Since each such member is divisible by 
$n \geq M^{3/4}\log^2{M}$, which is greater than the length of $I$, we must
have that all such members divide the same integer in $I$.  Thus, case A is 
true.\smallskip

To finish the proof of our lemma, we now show that 
$\text{\rm lcm\ }M(n) = P$ for all $n$ satisfying (15).
Fix an $n$ satisfying (15) and 
let $p^a \leq M^{1/4 - \epsilon}$ be the largest power of the prime
$p$ that is $\leq M^{1/4-\epsilon}$.  Let $p^e$  be the exact power
of $p$ which divides $n$.  Thus, $e \leq a$. 
We will show there exists an $m_j \in M(n)$ with 
$$
m_j = n p^{a-e}l_1l_2,\ \text{\rm where $l_1$ and $l_2$ are primes
with gcd$(l_1l_2,n) = 1$},
$$
which will imply that $m_j$ is divisible by $p^a$, and thus
$p^a | \text{\rm lcm }M(n)$.
Such an $m_j$ exists if we can just find primes 
$l_1, l_2 \leq M^{1/4-\epsilon}$ which satisfy 
$$
\sqrt{\frac{M}{np^{a-e}}} \leq l_1 < l_2 
\leq \sqrt{\left ( 1 + \frac{1}{\log{M}} \right )\frac{M}{np^{a-e}}},
\ \text{ gcd }(l_1l_2,n)=1.\eqno{(16)}
$$

To see that it is possible to find $l_1$ and $l_2$ we first observe that the
lower limit of the interval in (16) is
$$
\sqrt{\frac{M}{np^{a-e}}} \gg \sqrt{\frac{M}{(M^{3/4}\log^2{M})
 M^{1/4-\epsilon}}} = \frac{M^{\epsilon/2}}{\log{M}}, 
$$
and the length of the interval is the multiple $\sqrt{1 + \frac{1}{\log{M}}}
-1 \gg \frac{1} {\log{M}}$ of this lower limit.  By the Prime Number Theorem,
there are $\gg \frac{M^{\epsilon/2}}{\epsilon\log^3{M}}$ primes in this
interval,
and so for $M$ sufficiently large there must be two of them $l_1 < l_2$ which
do not divide $n < 2M^{3/4}\log^2{M}$.  These two
primes therefore satisfy (16).  To see that $l_1,l_2 < M^{1/4-\epsilon}$,
we observe that the upper limit of the interval in (16) satisfies
$$
\sqrt{\left ( 1 + {1 \over \log{M}} \right )
\frac{M}{np^{a-e}}} < \sqrt{\frac{2M}{n}} \leq
\sqrt{\frac{2M}{M^{3/4}\log^2{M}}} = \frac{\sqrt{2}M^{1/8}}{\log{M}}
< M^{1/4-\epsilon},
$$
for $M$ sufficiently large and $0 < \epsilon < 1/8$.  Thus, we can find
$l_1$ and $l_2$ as claimed, and so our lemma is proved.

\heading{VI.  Proof of Main Theorem}\endheading
We give here only a slightly more formal version of the proof outlined in
the introduction.\smallskip

Suppose we are given a rational number $r>0$ and an integer $N>r$.  Let
$M$ be the least integer where
$$
r \leq \sum_{N \leq n \leq M} {1 \over n} \leq r + {1 \over M}.
$$
Using the fact that $\sum_{1 \leq n \leq x} {1 \over n} = \log{x} + \gamma 
+ O(1/x)$, it is easy to see that $M/N = e^{r+O(1/N)}$.\smallskip\noindent
Using Proposition $1$ with $\epsilon = 1/6$ we have that for $N$ sufficiently
large, there are integers $d_1,...,d_l$ with
$$
N \leq d_1 < d_2 < \cdots < d_l < M = e^{r+O(1/N)}N,
$$
such that if 
$$
{u \over v}\ :=\ \sum_{N \leq n \leq M \atop n \neq d_1,...,d_l} {1 \over n}
= r - (3r + o(1)) {\log\log{N} \over \log{N}},\ \text{gcd}(u,v)=1,
$$
then where all the prime power factors of $v$ are $\leq N^{1/4-1/6} = 
N^{1/12}$.  Let
$$
{a \over b} = r - {u \over v} = (3r + o(1)) {\log\log{N} \over \log{N}},
\ \text{gcd}(a,b)=1.
$$
We observe that once $N$ is large enough, all the prime power factors of $b$
will be $\leq N^{1/12}$.  Invoking Proposition 2 with $\epsilon = 1/6$
we have that there are integers $n_1,...,n_k$ with
$$
M \leq n_1 < \cdots < n_k \leq e^{c\cdot a/b}M,
$$
where $c$ is some constant, and such that
$$
{a \over b} = {1 \over n_1} + \cdots + {1 \over n_k},
$$
Thus, we have the representation for $r$:
$$
r = \left ( \sum_{N \leq n \leq M \atop n \neq d_1,...,d_l} { 1 \over n}
\right ) + {1 \over n_1} + {1 \over n_2} + \cdots + {1 \over n_k},
$$
where
$$
n_k \leq e^{c\cdot a/b}M = \left \{1 +  (3cr + o(1)){\log\log{N} \over \log{N}}
\right \} M = \left \{ e^r + O_r\left ( {\log\log{N} \over \log{N}} 
\right ) \right \} N.
$$
This proves the first part of the Main Theorem.\smallskip

To see that the $O_r\left ( {\log\log{N} \over \log{N}}
\right )$ error term is best-possible, suppose that
$$
r = {a \over b} = {1 \over x_1} + \cdots + {1 \over x_k},\  \text{gcd$(a,b)=1$},
$$
where $N \leq x_1,...,x_k \leq cN$ are distinct integers, and let
$x$ be the largest of the $x_i$'s.  
We claim that the largest prime $p$ dividing the $x_i$'s satisfies
$p < {x \over \log{x}}(1+o(1))$.
To see this, let 
$$
x_1 = p m_1 < x_2 = p m_2 < \cdots < x_l = p m_l
$$  
be all the $x_i$'s divisible by $p$.  If $p | b$ then since $b$ remains
bounded as $x$ varies, we would have that $p \leq b < x/\log{x}$ once $x$
is large enough.  If, on the other hand, $p \nmid b$, then we must have
that $p \nmid b'$ either, where $b'$ is given by
$$
{a' \over b'} = {1 \over x_1} + \cdots + {1 \over x_l}
= {1 \over p} \left ( {1 \over m_1} + \cdots + {1 \over m_l} 
\right ),\ \text{ gcd$(a',b')=1$}.
$$
Thus, $p$ divides
$$
\eqalign{
\text{lcm}\{m_1,...m_l\} \left \{ {1 \over m_1} + \cdots + {1 \over m_l}
\right \} &\leq \text{lcm}\{2,3,...,m_l\} \left \{ 1 + {1 \over 2} + {1 \over 3}
+ \cdots + {1 \over m_l} \right \}\cr
&= e^{m_l(1+o(1))},
}
$$
and so,
$$
x \geq pm_l > p\log{p}(1+o(1));
$$
or in other words,
$$
p < {x \over \log{x}}(1+o(1)).
$$ 

Making use of this bound on $p$ we have that
$$
r \leq \sum_{N \leq n \leq cN\atop p |n \Longrightarrow 
p < {cN\over \log(cN)}(1+o(1))} {1 \over n} = 
\left (\sum_{N \leq n \leq cN} {1 \over n}\right ) - \left (
\sum_{N \leq mp \leq cN \atop p > {cN \over \log{cN}}(1+o(1)) } {1 \over mp}
\right ) 
$$
Applying lemma 2 to this last pair of terms, together with the estimate
$\sum_{n \leq x} {1 \over n} = \log{x} + \gamma + O(1/x)$, we find that
$$
r \leq \log{c} - (\log{c} + o(1)){\log\log{N} \over \log{N}}.
$$
Solving for $c$ we find that 
$$
c \geq e^{r} \left ( 1 +  {(r+o(1))\log\log{N} \over \log{N}}\right ).
$$
\bigskip\noindent
Acknowledgements:  I would like to thank Drs. Andrew Granville and Carl
Pomerance for their comments and suggestions.  I would also like to thank
Greg Martin for the enlightening conversations I had with him by email and
in person. 

\enddocument